\documentclass[12pt,leqno]{amsart}
\usepackage{amssymb,verbatim,enumerate,ifthen,cite}
\usepackage[mathscr]{eucal}
\usepackage[utf8]{inputenc}
\usepackage[T1]{fontenc}
\def\N{\mathbb{N}}
\def\R{\mathbb{R}}
\def\Q{\mathbb{Q}}

\newtheorem{theorem}{Theorem}[section]
\newtheorem*{theorem*}{Theorem}
\def\Thm#1#2{\ifthenelse{\equal{#1}{*}}{\begin{theorem*}#2\end{theorem*}}
             {\begin{theorem}\label{T#1}#2\end{theorem}}}

\def\thm#1{Theorem~\ref{T#1}}

\newtheorem{proposition}[theorem]{Proposition}
\newtheorem*{proposition*}{Proposition}
\def\Prp#1#2{\ifthenelse{\equal{#1}{*}}{\begin{proposition*}#2\end{proposition*}}
{\begin{proposition}\label{P#1}#2\end{proposition}}}

\newtheorem{corollary}[theorem]{Corollary}
\newtheorem*{corollary*}{Corollary}
\def\Cor#1#2{\ifthenelse{\equal{#1}{*}}{\begin{corollary*}#2\end{corollary*}}
             {\begin{corollary}\label{C#1}#2\end{corollary}}}

\newtheorem{lemma}[theorem]{Lemma}
\newtheorem*{lemma*}{Lemma}
\def\Lem#1#2{\ifthenelse{\equal{#1}{*}}{\begin{lemma*}#2\end{lemma*}}
             {\begin{lemma}\label{L#1}#2\end{lemma}}}
\def\lem#1{Lemma~\ref{L#1}}

\theoremstyle{definition}
\newtheorem{remark}[theorem]{Remark}
\newtheorem*{remark*}{Remark}
\def\Rem#1#2{\ifthenelse{\equal{#1}{*}}{\begin{remark}\rm #2\end{remark}}
             {\begin{remark}\label{R#1}\rm #2\end{remark}}}

\newtheorem{example}[theorem]{Example}
\newtheorem*{example*}{Example}
\def\Exa#1#2{\ifthenelse{\equal{#1}{*}}{\begin{example*}\rm #2\end{example*}}
             {\begin{example}\label{Ex#1}\rm #2\end{example}}}

\def\eq#1{{\rm(\ref{E#1})}}
\def\Eq#1#2{\ifthenelse{\equal{#1}{*}}
  {\begin{equation*}\begin{aligned}[]#2\end{aligned}\end{equation*}}
  {\begin{equation}\begin{aligned}[]\label{E#1}#2\end{aligned}\end{equation}}}

\begin{document}
\begin{flushright}
%\textit{Submitted to:} 
\end{flushright}
\vspace{5mm}

\date{\today}

\title{Decomposition of higher-order Wright convex functions revisited}

\author[Zs. P\'ales]{Zsolt P\'ales}
\address[Zs. P\'ales]{Institute of Mathematics, University of Debrecen, 
H-4002 Debrecen, Pf.\ 400, Hungary}
\email{pales@science.unideb.hu}
\author[M. K. Shihab]{Mahmood Kamil Shihab}
\address[M. K. Shihab]{Doctoral School of Mathematical and Computational Sciences, University of Debrecen, H-4002 Debrecen, Pf.\ 400, Hungary}
\email{mahmood.kamil@science.unideb.hu}


%\thanks{}

\begin{abstract}
In 2009, Maksa and P\'ales established an extension of the decomposition theorem of Ng in the context of higher-order convexity notions. They proved that a real function is Wright convex of order $n$ if and only if it can be decomposed as the sum of a convex function of order $n$ and a polynomial function of order at most $n$. Their proof was based on transfinite tools in the background. The main purpose of this paper is to adopt the methods of a paper of P\'ales published in 2020 and establish a new and elementary proof for the theorem of Maksa and P\'ales.
\end{abstract}

\maketitle

\section{Introduction} 

Throughout this paper, let $I$ denote a proper open subinterval of the real line. 
For a function $f:I\to\R$, one can define three notions of convexity in the following way:

\noindent -- $f$ is called \emph{Jensen convex} on $I$, if
\Eq{*}{
  f\bigg(\frac{x+y}{2}\bigg)\leq \frac{f(x)+f(y)}{2} \qquad (x,y\in I).
}
-- $f$ is called \emph{Wright convex} on $I$, if
\Eq{*}{
  f(tx+(1-t)y)+f((1-t)x+ty)\leq f(x)+f(y)\qquad (x,y\in I,\,t\in[0,1]).
}
-- $f$ is called \emph{convex} on $I$, if
\Eq{*}{
  f(tx+(1-t)y)\leq tf(x)+(1-t)f(y)\qquad (x,y\in I,\,t\in[0,1]).
}
One can see that convexity implies Wright convexity and Wright convexity implies Jensen convexity and non of these implications can be reversed. 
For an overview about the generalizations, stability and regularity properties of Wright convex functions, we refer to the list of references, which may give the impression that this subfield of functional equations and inequalities is still in the focus of research.

A celebrated result of C.\ T.\ Ng \cite{Ng87b} established a deeper connection between convexity and Wright convexity. It characterizes Wright convex functions as those functions that are of the form $f=g+a$, where $g$ is convex and $a$ is additive. The original proof of the paper \cite{Ng87b} applied de Bruijn's theorem \cite{Bru51} which is related to functions which have continuous differences. Several subsequent proofs of the result of Ng (c.f., Nikodem \cite{Nik89a} and Kominek \cite{Kom03}) used another approach, which was based on Rode's theorem \cite{Rod78}. Basically, all the previously known proofs used transfinite induction for the construction of the additive part $a$ of the decomposition. In a recent paper \cite{Pal20a},  P\'ales obtained a new proof in which the convex summand $g$ was first constructed in an elementary way. Therefore, there was no transfinite induction involved.

In the paper \cite{MakPal09b}, Maksa and P\'ales extended the decomposition theorem of Ng to the context of higher-order convexity notions. They proved that a real function is Wright convex of order $n$ if and only if it can be decomposed as the sum of a convex function of order $n$ and a polynomial function of order at most $n$. Their proof was again using transfinite tools in the background. The main purpose of this paper is to adopt the methods of the paper \cite{Pal20a} and establish a new and elementary proof for the theorem of Maksa and P\'ales.

\section{Higher-order convex and Wright convex functions}

In what follows, we are going to define several higher-order convexity concepts in terms of {\em difference operators} and {\em divided differences}. 

We recall that, for a fixed real number $h$, the operator $\Delta_h$, acting on a real function $f:I\to\R$, is defined by
\Eq{*}{
  \Delta_{h}f(x):=f(x+h)-f(x) \qquad(x\in I\cap(I-h)).
}
Obviously, if $|h|$ is small enough, then $I\cap(I-h)$ 
is a non-void open interval again. The product of these operators 
can also be defined in the usual way (see e.g.\ Kuczma \cite{Kuc85}).

Given a fixed $n\in\N$, a map $f:I\to \R$ is said to be {\em Jensen convex of order $n$} (briefly $n$-Jensen convex) if
\Eq{nJ}{
\Delta^{n+1}_{h}f(x)\geq 0 \qquad  \big(h>0,\, x\in I\cap(I-(n+1)h)\big).
}
A map $f:I\to\R$ is said to be {\em Wright convex of order $n$} (briefly $n$-Wright convex) if it satisfies the functional inequality
\Eq{nW}{
\Delta_{h_{1}}\cdots\Delta_{h_{n+1}} f(x)\geq 0 \qquad\big(h_1,\dots, h_{n+1}>0,\, x\in I\cap(I-(h_1+\dots+ h_{n+1}))\big).
}

In the investigation of functional inequalities \eq{nJ} and \eq{nW}, those maps that fulfil these inequalities with equality play a fundamental role in the theory of linear functional equations. Therefore, for $n\in\N$, we consider the equation
\Eq{*}{
\Delta^{n+1}_{h}f(x)=0 \qquad  (h>0,\,x\in I\cap(I-(n+1)h)),
}
which is termed the \emph{Fr\'echet functional equation} in this theory.
It is well-known (see \cite{Kuc85}, \cite{Sze91}) that $f:I\to\R$ satisfies this equation if and only if it is a \emph{polynomial function of degree at most $n$}, i.e.,
it has the representation
\Eq{*}{
f(x)=a_{0}+a_{1}(x)+\dots +a_{n}(x) \qquad  (x \in I),
}
where $a_{0}\in\R$ and $a_{k}$ is the {\em diagonalization} of some 
$k$-additive and symmetric function $A_{k}:\R^{k}\to\R$, that is, 
$a_{k}(x)=A_{k}(x,\dots, x), \,\,(x\in\R, \,\, k=1,\dots, n)$. 
Standard polynomials are exactly the continuous polynomial functions.
On the other hand, using Hamel bases, it is not difficult to construct non-continuous polynomial functions (see \cite{Kuc85}). 

The {\em divided difference} of  the function $f:I\to\R$ with respect to the
pairwise distinct points $x_0,\ldots,x_n\in I$ is defined by
\Eq{*}{
  [x_0,\ldots,x_n;f]=
  \sum_{i=0}^n\frac{f(x_i)}
      {\prod_{\begin{smallmatrix} j=0 \\ j\neq i \end{smallmatrix}}^n(x_i-x_j)}.
}
Obviously, divided differences are symmetric functions of their variables,
furthermore, it is easy to show that they enjoy the following recursive property
\Eq{*}{
  [x_0,\ldots,x_n;f]=\frac{[x_1,\ldots,x_n;f]-[x_0,\ldots,x_{n-1};f]}{x_n-x_0}
}
for all $n\in\N$ and pairwise distinct elements $x_0,\ldots,x_n\in I$.

Based on the works of T.~Popoviciu \cite{Pop34,Pop44}, given $n\in\N$, a map $f:I\to\R$ is said to be \emph{convex of order $n$ on $I$} (shortly $n$-convex on $I$) if the inequality
\Eq{S2.2}{
  [x_0,x_1,\ldots,x_n,x_{n+1};f]\geq0
}
holds for all pairwise distinct elements $x_0,x_1,\ldots,x_n,x_{n+1}\in I$.
Due to the symmetry of divided differences, without loss of generality, we may assume $x_0<x_1<\cdots<x_n<x_{n+1}$ here.

The following result was obtained in the book \cite{Kuc85} and in a more general form in the paper \cite{GilPal08}.

\Lem{L1}{Let $n\in\N$. Then every $n$-convex function is $n$-\!Wright convex, and every $n$-\!Wright convex function is $n$-Jensen convex.}

One of the main results of the paper \cite{MakPal09b} established the following generalization of Ng's decomposition theorem \cite{Ng87b}.

\Thm{T1}{Let $n\in\N$ and $f:I\to\R$ be an $n$-\!Wright convex function. Then, there exist an $n$-convex function $g:I\to\R$ and a polynomial function $P:\R\to\R$ of degree at most $n$ such that
\Eq{*}{
  f(x)=g(x)+P(x)\qquad(x\in I).
}}

Our aim is to obtain a new proof for this result.

\Lem{L2}{
	Let $n\in\N$ and $f:I\to\R$ be an $n$-Jensen convex function. Then there exists a continuous $n$-convex function $g:I\to\R$ such that $g\vert_{I\cap\Q}=f\vert_{I\cap\Q}$.
}
\begin{proof}
    Using \cite[Lemma 5.1]{GilPal08}, we have the following identity 
    \Eq{*}{
      [x,x+h,\dots,x+nh,x+(n+1)h;f]=\frac{\Delta_h^{n+1}f(x)}{(n+1)!h^{n+1}} \qquad(h>0,\,x\in I\cap(I-(n+1)h)).
    }
	Therefore, the $n$-Jensen convexity of $f$ implies that 
	\Eq{*}{
	  [x,x+h,x+2h\dots,x+(n+1)h;f]\geq 0 
	  \qquad(h>0,\,x\in I\cap(I-(n+1)h)).
	}
	In the terminology of the paper \cite{GilPal08}, this property says that $f$ is $(t_0,\dots,t_n)$-convex with $t_0 =\dots=t_n=1$. According to \cite[Theorem 3.2]{GilPal08}, it follows that $f$ is $(r_0,\dots,r_n)$-convex for all positive rational numbers $(r_0,\dots,r_n)$-convex, that is
	\Eq{rc1}{
	  [x,x+r_0h,x+(r_0+r_1)h,\dots,x+(r_0+\dots+r_n)h;f]\geq 0 \qquad\qquad \\
	  (h>0,\,x\in I\cap(I-(r_0+\dots+r_n)h)).
	}
	We now deduce that $f$ satisfies the $n$-convexity property with rational arguments, that is,
	\Eq{rc2}{
	  [x_0,x_1,\dots,x_n,x_{n+1};f]\geq 0 
	  \qquad(x_0,\dots,x_{n+1}\in I\cap\Q \mbox{ such that $x_i\neq x_j$ if $i\neq j$}  ).
	}
	Indeed, let $x_0,x_1,\dots,x_n,x_{n+1}\in I\cap\Q$ be arbitrary. Without loss of generality, we may assume that $x_0<x_1<\dots<x_n<x_{n+1}$. Applying now the inequality \eq{rc1} with $x:=x_0$, $h=1$ and $r_i:=x_{i+1}-x_i$ for $i\in\{0,\dots,n\}$, we can see that \eq{rc2} holds.
	    
    Claim: For any compact subinterval $[a,b]\subseteq I$, there exists $L\geq 0$ such that
    \Eq{*}{
      |f(x)-f(y)|\leq L|x-y| \qquad(x,y\in[a,b]\cap\Q).
    }
    
    To show this claim, let $[a,b]\subseteq I$. Without loss of generality, we can assume that $a,b\in\Q$. Let $x,y\in[a,b]\cap\Q$ with $x<y$ be arbitrary. Then, for all pairwise distinct elements $u_1,\dots,u_n\in (I\setminus[a,b])\cap\Q$, we get
    \Eq{ddi}{
	  [y,x,u_1,\dots,u_n;f]\geq 0,
	}
	that is,
	\Eq{*}{
	  \frac{1}{y-x} \bigg(\frac{f(y)}{\prod_{j=1}^n(y-u_j)}
	  -\frac{f(x)}{\prod_{j=1}^n(x-u_j)}\bigg)
	  \geq -\sum_{i=1}^n\frac{f(u_i)}{(u_i-x)(u_i-y)
	        \prod_{j\in\{1,\dots,n\}\setminus\{i\}}(u_i-u_j)}.
	}
	The mapping 
	\Eq{*}{
	[a,b]^2\ni(x,y)\mapsto -\sum_{i=1}^n\frac{f(u_i)}{(u_i-x)(u_i-y)
	        \prod_{j\in\{1,\dots,n\}\setminus\{i\}}(u_i-u_j)}
	}
	is continuous on $[a,b]^2$ and therefore it is bounded from below by a
	constant $C(u_1,\dots,u_n)$, therefore, for all $x,y\in[a,b]\cap\Q$, the inequality \eq{ddi} implies 
    \Eq{ineq}{
	  \frac{1}{y-x} \bigg(\frac{f(y)}{\prod_{j=1}^n(y-u_j)}
	  -\frac{f(x)}{\prod_{j=1}^n(x-u_j)}\bigg)
	  \geq C(u_1,\dots,u_n).
	}
	Let first $u_1<\dots<u_n<a$ be fixed elements of $I\cap\Q$ and define $U(t):=\prod_{j=1}^n(t-u_j)$ for $t\in[a,b]$. Then $U$ is an increasing and positive polynomial on $[a,b]$. Therefore $U\leq U(b)$ and there exists a positive number $M>0$ such that $|U'|\leq M$ on $[a,b]$. Hence, by the Lagrange mean value theorem, $U$ is Lipschitzian over $[a,b]$ with a Lipschitz modulus $M$. From \eq{ineq}, for all $x,y\in[a,b]\cap\Q$, it follows that
    \Eq{xy}{
	  \frac{1}{y-x} \bigg(\frac{f(y)}{U(y)}
	  -\frac{f(x)}{U(x)}\bigg)
	  \geq C(u_1,\dots,u_n).
    }
    By putting $x=a$, this inequality yields
    \Eq{*}{
    f(y)&\geq \frac{U(y)}{U(a)}f(a)+ U(y)C(u_1,\dots,u_n)(y-a)\\
      &\geq -\frac{U(y)}{U(a)}|f(a)|-U(y)|C(u_1,\dots,u_n)|(y-a)\\
      &\geq -\frac{U(b)}{U(a)}|f(a)|-U(b)|C(u_1,\dots,u_n)|(b-a),
    }
    which shows that $f$ is bounded from below on $[a,b]\cap\Q$. On the other hand, putting $y=b$ in \eq{xy}, we can obtain that
    \Eq{*}{
      f(x)&\leq \frac{U(x)}{U(b)}f(b)+ U(x)C(u_1,\dots,u_n)(x-b)\\
      &\leq \frac{U(x)}{U(b)}|f(b)|+ U(x)|C(u_1,\dots,u_n)|(b-x)\\
      &\leq |f(b)|+ U(b)|C(u_1,\dots,u_n)|(b-a),
    }
    which shows that $f$ is also bounded from above on $[a,b]\cap\Q$. Thus, there exists a nonnegative number $K$ such that $|f(x)|\leq K$ for $x\in[a,b]\cap\Q$.
    The inequality \eq{xy} now yields
    \Eq{*}{
    f(y)-f(x)&\geq \frac{U(y)-U(x)}{U(x)}f(x)+ U(y)C(u_1,\dots,u_n)(y-x)\\
    &\geq -\frac{U(y)-U(x)}{U(x)}|f(x)|- U(y)|C(u_1,\dots,u_n)|(y-x)\\
    &\geq -\frac{M(y-x)}{U(a)}K- U(b)|C(u_1,\dots,u_n)|(y-x)\\
    &= -\bigg(\frac{MK}{U(a)}+U(b)|C(u_1,\dots,u_n)|\bigg)(y-x).
}

    Let, additionally $b<u_n'$. Then $V(t):=(t-u_n')\prod_{j=1}^{n-1}(t-u_j)$ for $t\in[a,b]$. Then $V$ is negative polynomial on $[a,b]$. Therefore there exist  positive numbers $M_0,M_1$ and $M_2$ such that $M_0\leq |V|\leq M_1$ and $|V'|\leq M_2$ on $[a,b]$. Hence again by the Lagrange mean value theorem $V$ is Lipschitzian over $[a,b]$ with a Lipschitz modulus $M_2$. From \eq{ineq} for all $x,y\in [a,b]\cap\Q$ we have
    \Eq{*}{
	  \frac{1}{y-x} \bigg(\frac{f(y)}{V(y)}
	  -\frac{f(x)}{V(x)}\bigg)
	  \geq C(u_1,\dots,u_{n-1},u_n'),
    }
    which is equivalent to
    \Eq{*}{
f(y)-f(x)&\leq \frac{V(y)-V(x)}{V(x)}f(x)+V(y)C(u_1,\dots,u_{n-1},u'_n)(y-x)\\
&\leq \bigg|\frac{V(y)-V(x)}{V(x)}f(x)\bigg|+|V(y)C(u_1,\dots,u_{n-1},u'_n)|(y-x)\\
&\leq \frac{M_2|y-x|}{M_0}K+M_1|C(u_1,\dots,u_{n-1},u'_n)|(y-x)\\
&=\bigg(\frac{M_2K}{M_0}+M_1|C(u_1.\dots,u_{n-1},u'_n)|\bigg)(y-x)
}
Define
\Eq{*}{
	L:=\max\bigg\{  \bigg(\frac{MK}{U(a)}+U(b)|C(u_1,\dots,u_n)|\bigg),\bigg(\frac{M_2K}{M_0}+M_1|C(u_1,\dots,u_{n-1},u'_n)|\bigg)\bigg\}.
}
Therefore $f$ is Lipschitz with modulus $L$ on the dense set $[a,b]\cap\Q$.
By applying \cite[Lemma 1]{Pal20a} with $D=[a,b]\cap\Q$, we get that there exists a continuous function $g:I\to \R$ such that $g|_{I\cap\Q}=f|_{I\cap\Q}$. 

To complete the proof, we have to show that $g$ is $n$-convex. Let $y_0,\dots,y_{n+1}$ be arbitrary pairwise distinct elements of $I$. Then, for all $j\in\{0,\dots,n+1\}$, there exist rational sequences $(x_{k,j})_{k\in\N}$ converging to $y_j$ as $k\to\infty$ with the property that the elements $x_{k,0},x_{k,1},\dots,x_{k,n},x_{k,n+1}$ are pairwise distinct for all $k\in\N$. Then, by applying the $n$-convexity property of $f$ with rational arguments (i.e., inequality \eq{rc2}), for all $k\in\N$, we get
\Eq{*}{
  [x_{k,0},x_{k,1},\dots,x_{k,n},x_{k,n+1};g]
  =[x_{k,0},x_{k,1},\dots,x_{k,n},x_{k,n+1};f]\geq0
}
By the continuity of $g$ the $(n+1)$st-order divided difference of $g$ is a continuous function of its arguments. Thus, upon taking the limit as $k\to\infty$, the above inequality yields that
\Eq{*}{
  [y_0,y_1,\dots,y_n,y_{n+1};g]\geq0.
}
Therefore, $g$ is $n$-convex, indeed.
\end{proof}

Now we give a new proof for \thm{T1}.
\begin{proof}[Proof of \thm{T1}]
	Since $f:I\to\R$ is $n$-Wright convex, therefore by \lem{L1}, $f$ is $n$-Jensen convex, i.e, $f$ satisfies \eq{nJ}. In particular, $f|_{I\cap\Q}$ is $n$-Jensen convex on $I\cap\Q$. Thus, in view of \lem{L2}, there exists a continuous $n$-convex function $g:I\to\R$ such that $f|_{I\cap\Q}=g|_{I\cap\Q}$.
	
	To complete the proof, we show that $f-g$ is a polynomial function of degree at most $n$. For this we prove that $\Delta_h^{n+1}(f-g)(x)=0$ for $h>0$ and $x\in I\cap(I-(n+1)h)$. This equation is equivalent to
	\Eq{eq7}{
		\Delta_h^{n+1}f(x)=\Delta_h^{n+1}g(x)\qquad (h>0, x\in I\cap(I-(n+1)h)).
	}
	More generally, we will show that
	\Eq{eq8}{
	  \Delta_{h_0}\cdots\Delta_{h_n}f(x)
	  =\Delta_{h_0}\cdots\Delta_{h_n}g(x)
	}
	holds for all $h_0,\dots,h_n>0$ and $x\in I\cap(I-(h_0+\dots+h_n))$.
	
	By the $n$-Wright convexity of $f$, for all $h_0,\dots,h_n>0$ and $x\in I\cap(I-(h_0+\dots+h_n))$, we have the inequality
	\Eq{eq9}{
		\Delta_{h_0}\cdots\Delta_{h_n}f(x)\geq 0.
	}
	This implies that $\Delta_{h_1}\cdots\Delta_{h_n}f:I\cap (I-(h_1+\dots+h_n))\to\R$ is nondecreasing for all $h_1,\dots,h_{n}>0$. On the other hand, $\Delta_{h_1}\cdots\Delta_{h_{n}}g:I\cap (I-(h_1+\dots+h_n))\to\R$ is continuous and the equality $f|_{I\cap\Q}=g|_{I\cap\Q}$ gives us $\Delta_{h_1}\cdots\Delta_{h_{n}}f(x)=\Delta_{h_1}\cdots\Delta_{h_{n}}g(x)$ for $x\in I\cap (I-(h_1+\dots+h_n))\cap\Q$ whenever $h_1,\dots,h_n\in\Q_+$. Applying \cite[Lemma 3]{Pal20a}, it follows that these two functions are equal to each other also at irrational points of $I\cap (I-(h_1+\dots+h_n))$, that is, 
	\Eq{*}{
		\Delta_{h_1}\cdots\Delta_{h_{n}}f(x)=\Delta_{h_1}\cdots\Delta_{h_{n}}g(x)\qquad (h_1,\dots,h_n\in\Q_+, x\in I\cap (I-(h_1+\dots+h_n)).
	}
	Applying this equality at $x+h_0$ and at $x$, and then subtracting the two equalities side by side, we can see that \eq{eq8} is valid if $h_0>0$, $h_1,\dots,h_n\in\Q_+$ and $x\in I\cap(I-(h_0+\dots+h_n))$. 
	
	Let $k\in\{0,\dots,n\}$ and consider the statement $S_k$ which says that \eq{eq8} holds if $h_0,\dots,h_k>0$, $h_{k+1},\dots,h_n\in\Q_+$ and $x\in I\cap(I-(h_0+\dots+h_n))$. According to the previous argument, we have that $S_0$ is true. Now assume that, for some $k\in\{0,\dots,n-1\}$, the statement $S_k$ holds. We show that $S_{k+1}$ is also valid. To prove this let $h_0,\dots,h_{k+1}>0$, $h_{k+2},\dots,h_n\in\Q_+$ and $x\in I\cap(I-(h_0+\dots+h_n))$ and let $h_{k+1}'<h_{k+1}$ be an arbitrary rational number. Then, by the $n$-Wright convexity of $f$, we have
	\Eq{*}{
	  \bigg(\prod_{i\in\{0,\dots,n\}\setminus\{k+1\}}\Delta_{h_i}\bigg)\Delta_{h_{k+1}-h'_{k+1}}f(x+h'_{k+1})\geq0.
	}
	Therefore, using the statement $S_k$ in the last step, we get
	\Eq{*}{
	  \Delta_{h_0}\cdots\Delta_{h_n}f(x)
	  &=\bigg(\prod_{i\in\{0,\dots,n\}\setminus\{k+1\}}\Delta_{h_i}\bigg)\cdot\Delta_{h_{k+1}}f(x)\\
	  &=\bigg(\prod_{i\in\{0,\dots,n\}\setminus\{k+1\}}\Delta_{h_i}\bigg)\big(f(x+h_{k+1})-f(x)\big)\\
	  &=\bigg(\prod_{i\in\{0,\dots,n\}\setminus\{k+1\}}\Delta_{h_i}\bigg)\big(f(x+h_{k+1})-f(x+h'_{k+1})+f(x+h'_{k+1})-f(x)\big)\\
	  &=\bigg(\prod_{i\in\{0,\dots,n\}\setminus\{k+1\}}\Delta_{h_i}\bigg)\big(\Delta_{h_{k+1}-h'_{k+1}}f(x+h'_{k+1})+\Delta_{h'_{k+1}}f(x)\big)\\
	  &\geq\bigg(\prod_{i\in\{0,\dots,n\}\setminus\{k+1\}}\Delta_{h_i}\bigg)\Delta_{h'_{k+1}}f(x)\\
	  &=\bigg(\prod_{i\in\{0,\dots,n\}\setminus\{k+1\}}\Delta_{h_i}\bigg)\Delta_{h'_{k+1}}g(x).
	}
	Using that $g$ is continuous, after taking the limit $h_{k+1}'\to h_{k+1}$, we get that
	\Eq{*}{
	  \Delta_{h_0}\cdots\Delta_{h_n}f(x)
    \geq \bigg(\prod_{i\in\{0,\dots,n\}\setminus\{k+1\}}\Delta_{h_i}\bigg)\Delta_{h_{k+1}}g(x)=\Delta_{h_0}\cdots\Delta_{h_n}g(x).
	}
	To prove the other direction let $h''_{k+1}>h_{k+1}$ be an arbitrary rational number, again by the $n$-Wright convexity of $f$ and statement $S_k$, we get
	\Eq{*}{
		\Delta_{h_0}\cdots\Delta_{h_n}f(x)
		&=\bigg(\prod_{i\in\{0,\dots,n\}\setminus\{k+1\}}\Delta_{h_i}\bigg)\cdot\Delta_{h_{k+1}}f(x)\\
		&=\bigg(\prod_{i\in\{0,\dots,n\}\setminus\{k+1\}}\Delta_{h_i}\bigg)\big(f(x+h_{k+1})-f(x)\big)\\
		&=\bigg(\prod_{i\in\{0,\dots,n\}\setminus\{k+1\}}\Delta_{h_i}\bigg)\big(f(x+h_{k+1})-f(x+h''_{k+1})+f(x+h''_{k+1})-f(x)\big)\\
		&=\bigg(\prod_{i\in\{0,\dots,n\}\setminus\{k+1\}}\Delta_{h_i}\bigg)\big(-\Delta_{h''_{k+1}-h_{k+1}}f(x+h_{k+1})+\Delta_{h''_{k+1}}f(x)\big)\\
		&\leq\bigg(\prod_{i\in\{0,\dots,n\}\setminus\{k+1\}}\Delta_{h_i}\bigg)\Delta_{h''_{k+1}}f(x)\\
		&=\bigg(\prod_{i\in\{0,\dots,n\}\setminus\{k+1\}}\Delta_{h_i}\bigg)\Delta_{h''_{k+1}}g(x).
	}
	Upon taking the limit $h_{k+1}''\to h_{k+1}$ and using that $g$ is continuous, we get that
	\Eq{*}{
		\Delta_{h_0}\cdots\Delta_{h_n}f(x)
		\leq \bigg(\prod_{i\in\{0,\dots,n\}\setminus\{k+1\}}\Delta_{h_i}\bigg)\Delta_{h_{k+1}}g(x)=\Delta_{h_0}\cdots\Delta_{h_n}g(x).
	}
Combining the above two inequalities, we can see that $S_{k+1}$ is valid and hence, we have proved that $S_k$ is true for all $k\in\{0,\dots,n\}$, in particular, for $k=n$. This means that \eq{eq8} is satisfied for any $h_0,\cdots,h_n>0$. Consequently, it holds with $h_0=\cdots=h_n=h>0$, proving that \eq{eq7} is satisfied. This implies that $f-g$ is a polynomial function of degree at most $n$.
\end{proof}

\section{Statements and Declarations}

\textbf{Conflict of Interests:} None

\textbf{Funding Information:} The research of the first author was supported by the K-134191 NKFIH Grant and the 2019-2.1.11-TÉT-2019-00049, the EFOP-3.6.1-16-2016-00022 and the EFOP-3.6.2-16-2017-00015 projects. The last two projects are co-financed by the European Union and the European Social Fund.

\textbf{Code availability:} Not applicable

\textbf{Data availability:} Not applicable

%\nocite{Pal20a}
%\bibliography{convexity}
%\bibliography{publ,funcequ}
%\bibliographystyle{amsplain}
%\end{document}
\def\MR#1{}

\end{document}